\documentclass[11pt, reqno]{amsart}
\usepackage{graphics}
\usepackage{amsfonts}

\def\Box{\vcenter{\vbox{\hrule\hbox{\vrule
     \vbox to 8.8pt{\hbox to 10pt{}\vfill}\vrule}\hrule}}}

\def\mymedskip{\vskip\medskipamount}
\def\mymedbreak{\par \ifdim\lastskip<\medskipamount
  \removelastskip \penalty-100 \mymedskip \fi}
\def\myaftermedspace{\par \ifdim\lastskip<\medskipamount
  \removelastskip \penalty55\mymedskip\fi}
\newcommand{\eop}{{\unskip\nobreak\hfil\penalty50
          \hskip2em\hbox{}\nobreak\hfil$\Box$
          \parfillskip=0pt \finalhyphendemerits=0 \par}}
\renewenvironment{proof}%
{\mymedbreak{\noindent\bf Proof:\enspace}}{\eop\myaftermedspace}
\newenvironment{proofofteor}[1]%
{\mymedbreak{\noindent\bf Proof of Theorem #1:\enspace}}{\eop\myaftermedspace}
\newtheorem{teor}{Theorem}[section]

\newtheorem{lem}[teor]{Lemma}
\newtheorem{cor}[teor]{Corollary}
\newtheorem{con}[teor]{Conjecture}

\newtheorem{rem}[teor]{Remark}
\newcommand{\beq}{\begin{equation}}
\newcommand{\eeq}{\end{equation}}
\newcommand{\beql}[1]{\begin{equation} \label{#1}}
\newcommand{\eeql}{\end{equation}}
\newcommand{\beqa}{\begin{eqnarray*}}
\newcommand{\eeqa}{\end{eqnarray*}}
\newcommand{\beqal}[1]{\begin{eqnarray} \label{#1}}
\newcommand{\eeqal}{\end{eqnarray}}
\newcommand{\beqan}{\begin{eqnarray}}
\newcommand{\eeqan}{\end{eqnarray}}
\newcommand{\bpf}{\begin{proof}}
\newcommand{\epf}{\end{proof}}
\newcommand{\bpft}[1]{\begin{proofofteor}{#1}}
\newcommand{\epft}{\end{proofofteor}}

\newcommand{\Ff}{{\mathbb F}}
\newcommand{\Zz}{{\mathbb Z}}
\newcommand{\Cc}{{\mathbb C}}

\newcommand{\gc}{\gamma}

\def\Tr{\operatorname{Tr}}

\def\Tr{\operatorname{Tr}}

\begin{document}

\title{Proofs of two conjectures on ternary weakly regular bent functions}

\author[Helleseth, Hollmann, Kholosha, Wang, and Xiang]
{Tor Helleseth$^1$, Henk D. L. Hollmann, Alexander Kholosha$^1$,\\
Zeying Wang, and Qing Xiang$^2$}
\thanks{$^1$Research supported by the Norwegian Research Council.}
\thanks{$^2$Research supported in part by NSF Grant DMS 0701049.}


\date{}

\begin{abstract}
We study ternary monomial functions of the form $f(x)=\Tr_n(ax^d)$,
where $x\in \Ff_{3^n}$ and $\Tr_n: \Ff_{3^n}\rightarrow \Ff_3$ is
the absolute trace function. Using a lemma of Hou \cite{hou},
Stickelberger's theorem on Gauss sums, and certain ternary weight
inequalities, we show that certain ternary monomial functions
arising from \cite{hk1} are weakly regular bent, settling a
conjecture of Helleseth and Kholosha \cite{hk1}. We also prove that
the Coulter-Matthews bent functions are weakly regular.
\end{abstract}

\maketitle

\section{Introduction and Summary of results}\label{intro}

Let $p$ be a prime, $n\geq 1$ be an integer. We will use $\Ff_{p^n}$
to denote the finite field of size $p^n$, and $\Ff_{p^n}^*$ to
denote the set of nonzero elements of $\Ff_{p^n}$. Let $f:
\Ff_{p^n}\rightarrow \Ff_p$ be a function. The {\it Walsh (or
Fourier) coefficient} of $f$ at $b\in \Ff_{p^n}$ is defined by
\begin{eqnarray*}
S_f(b)&=&\sum_{x \in\Ff_{p^n}}\omega^{f(x)-\Tr_n(bx)}
\end{eqnarray*}
where $\Tr_n: \Ff_{p^n}\rightarrow \Ff_p$ is the absolute trace
function, $\omega=e^{\frac{2\pi i}{p}}$ is a primitive complex $p$th
root of unity, and elements of $\Ff_{p}$ are considered as integers
modulo $p$. In the sequel, $S_a(b)$ is also used to denote the Walsh
transform coefficient of a function that depends on parameter $a$
when it is clear from the context which function we mean. The
function $f$ is said to be a {\it $p$-ary bent function} (or a {\it
generalized bent function}) if all its Walsh coefficients satisfy
$$|S_f(b)|^2=p^n.$$
A $p$-ary bent function $f$ is said to be {\it regular} if for every
$b \in\Ff_{p^n}$ the normalized Walsh coefficient
$p^{-\frac{n}{2}}S_f(b)$ is equal to a complex $p$th root of unity,
i.e., $p^{-\frac{n}{2}}S_f(b)=\omega^{f^*(b)}$ for some function
$f^*:\Ff_{p^n}\rightarrow \Ff_{p}$. A bent function $f$ is said to
be {\it weakly regular} if there exists a complex number $u$ with
$|u|=1$ such that $up^{-\frac{n}{2}}S_f(b)=\omega^{f^*(b)}$ for all
$b \in\Ff_{p^n}$, where $f^*:\Ff_{p^n}\rightarrow \Ff_{p}$ is a
function. In such a situation, the function $f^*$ is also a weakly
regular bent function and it is called the {\it dual} of $f$.

Binary bent functions are usually called {\it Boolean bent
functions}, or simply {\it bent functions}. These functions were
first introduced by Rothaus \cite{roth} in 1976. Later Kumar,
Scholtz, and Welch \cite{ksw} generalized the notion of a Boolean
bent function to that of a $p$-ary bent function. All known $p$-ary
bent functions but one are weakly regular. The only known example of bent
but not weakly regular bent function was constructed by Helleseth
and Kholosha \cite{hk1}.

Bent functions, and in general, $p$-ary bent functions are closely
related to other combinatorial and algebraic objects such as
Hadamard difference sets in $(\Ff_{2^n},+)$ \cite{dillon}, relative
difference sets \cite{pott}, planar functions, and commutative
semifields \cite{coulter, carletding, wqwx}. For future use, we explicitly state
the relationship between planar functions and $p$-ary bent functions
here. A function $F:\Ff_{p^n}\rightarrow\Ff_{p^n}$ is said to be
{\it planar} if the function from $\Ff_{p^n}$ to $\Ff_{p^n}$ induced
by the polynomial $F(X+a)-F(X)-F(a)$ is bijective for every nonzero
$a\in\Ff_{p^n}$. The following lemma gives the relationship between
planar functions and $p$-ary bent functions.

\begin{lem}\label{planarbent}
Let $F:\Ff_{p^n}\rightarrow\Ff_{p^n}$ be a function. Then $F$ is
planar if and only if $\Tr_n(aF(x))$ is $p$-ary bent for all $a\in
\Ff_{p^n}^*$.
\end{lem}

The proof of the lemma is fairly straightforward, see
\cite{carletd}. Almost all known planar functions
$F:\Ff_{p^n}\rightarrow\Ff_{p^n}$ are of Dembowski-Ostrom type,
namely the corresponding polynomials $F(X)$ have the form
$F(X)=\sum_{i,j}a_{ij}X^{p^i+p^j}\in \Ff_{p^n}[X]$. The
Coulter-Matthews planar functions are special since they are not of
Dembowski-Ostrom type. These planar functions can be defined as
follows. Let $n,k\geq 1$ be integers such that $\gcd(k,n)=1$. Then
the function $F: \Ff_{3^n}\rightarrow \Ff_{3^n}$ defined by
$$F(x)=x^{\frac {3^k +1} {2}},\;\forall x\in \Ff_{3^n}$$
is planar. Thus by Lemma~\ref{planarbent}, $\Tr_n(ax^{\frac {3^k +1}
{2}})$ is $3$-ary bent for every nonzero $a\in\Ff_{3^n}$. These bent
functions are usually called {\it the Coulter-Matthews bent
functions}. It is conjectured that the Coulter-Matthews bent
functions are weakly regular \cite{hk1}, \cite{newhou}. (Strictly
speaking, it was only stated as an open problem in \cite{hk1} to
decide whether the Coulter-Matthews bent functions are weakly
regular or not. But most people believed that these functions are
weakly regular bent.) In a recent paper \cite{newhou}, it was proved
that the Coulter-Matthews bent functions are weakly regular in two
special cases. We confirm the conjecture in this paper. Therefore
our first result in this paper is

\begin{teor}\label{minor}
Let $n,k\geq 1$ be integers such that $\gcd(k,n)=1$. Then the bent
function $\Tr_n(ax^{\frac {3^k +1} {2}})$, $a\in \Ff_{3^n}^*$, is
weakly regular bent.
\end{teor}

Helleseth and Kholosha \cite{hk1} surveyed all proven and
conjectured classes of $p$-ary monomial bent functions
$f:\Ff_{p^n}\rightarrow \Ff_p$ of the form $f(x)=\Tr_n(ax^d)$, where
$a\in\Ff_{p^n}^*$, $d$ is an integer, and $p$ is odd. (See Table 1
in \cite{hk1}.) In that paper, besides mentioning that it is an open
problem to decide whether the Coulter-Matthews bent functions are
weakly regular, the authors also made the following conjecture.

\begin{con} {\rm (\cite{hk1})}\label{hk}
Let $n=2k$ with $k$ odd. Then the ternary function $f$ mapping
$\Ff_{3^n}$ to $\Ff_3$ and given by
\[f(x)={\rm
Tr}_n\Big(ax^{\frac{3^n-1}{4}+3^k+1}\Big)\] is a weakly regular bent
function if $a=\xi^{\frac{3^k+1}{4}}$ and $\xi$ is a primitive
element of $\Ff_{3^n}$. Moreover, for $b\in\Ff_{3^n}$ the
corresponding Walsh transform coefficient of $f(x)$ is equal to
\[S_f(b)=-3^k\omega^{\pm{\rm
Tr}_k\left(\frac{b^{3^k+1}}{a(I+1)}\right)}\] where $I$ is a
primitive fourth root of unity in $\Ff_{3^n}$.
\end{con}

We will show that the ternary functions in the above conjecture are
indeed weakly regular bent. It still remains to prove the second
part of the conjecture. We state our second result in this paper as

\begin{teor}\label{major}
Let $k$ be an odd positive integer, and let $n=2k$. Then the ternary
function $f:\Ff_{3^n}\rightarrow \Ff_3$ defined by
$$f(x)=\Tr_n(ax^{\frac{3^n-1}{4}+3^k+1}),\;\forall x\in \Ff_{3^n},$$
is a weakly regular bent function if $a=\xi^{\frac{3^k+1}{4}}$ and
$\xi$ is a primitive element of $\Ff_{3^n}$.
\end{teor}

Our proofs of Theorem~\ref{minor} and \ref{major} rely on a lemma of
Hou \cite{hou}. The idea is of a $p$-adic nature, and it has been
used successfully a few times in the literature (see for example,
\cite{hx}, \cite{dwx}): Given a function $f:\Ff_{p^n}\rightarrow
\Ff_p$, it is usually difficult to compute the Walsh coefficients
$S_f(b)$ explicitly; sometimes, even computing the absolute values
of $S_f(b)$ is difficult. However, such difficulties can sometimes
be bypassed by divisibility considerations. To this end, we first
introduce Gauss sums, Stickelberger's theorem on Gauss sums, and
Hou's lemma.

\section{The Teichm\"uller character, Gauss sums, Stickelberger's Theorem,
and Hou's lemma}\label{character}

Let $p$ be a prime, $q=p^n$, and $n\geq1$. Let $\omega=e^{\frac{2\pi
i}{p}}$ be a primitive complex $p$th root of unity and let $\Tr_n$
be the trace from $\Ff_q$ to $\Zz/p\Zz$. Define
$$ \psi: \Ff_q \rightarrow \Cc^*, \quad \psi(x)=\omega^{\Tr_n(x)},$$
which is easily seen to be a nontrivial character of the additive
group of $\Ff_q$. Let
$$\chi:\Ff_q^* \rightarrow \Cc^* $$
be a character of $\Ff_q^*$ (the cyclic multiplicative group of
$\Ff_q$). We define the {\it Gauss sum} by
$$ g(\chi)=\sum_{a \in \Ff_q^*} \chi(a)\psi(a).$$
Note that if $\chi_0$ is the trivial multiplicative character of
$\Ff_q$, then $g(\chi_0)=-1$. Gauss sums can be viewed as the
Fourier coefficients in the Fourier expansion of $\psi|_{\Ff_q^*}$
in terms of the multiplicative characters of
 $\Ff_q$. That is, for every $c\in \Ff_q^*$,
\begin{equation}\label{inv}
\psi(c)=\frac{1} {q-1}\sum_{\chi\in X}g(\chi)\chi^{-1}(c),
\end{equation}
where $X$ denotes the character group of $\Ff_q^*$.

One of the elementary properties of Gauss sums is
\cite[Theorem~1.1.4]{be1}
\begin{equation}\label{eq2.2}
g(\chi)\overline{g(\chi)}=q,\hspace{0.1in}{\rm if}\hspace{0.1in}
\chi\neq \chi_0.
\end{equation}

A deeper result on Gauss sums is Stickelberger's theorem
(Theorem~\ref{stick} below) on the prime ideal factorization of
Gauss sums. We first introduce some notation. Let $a$ be any integer
not divisible by $q-1$. Then there are {\em unique\/} integers $a_0,
\ldots, a_{n-1}$ with $0\leq a_i\leq p-1$ for all $i$, $0\leq i\leq
n-1$ such that
\[a\equiv a_0+a_1p+\cdots +a_{n-1}p^{n-1} (\bmod q-1).\]
We define the ($p$-ary) {\em weight} of $a$ (mod $q-1$), denoted by $w(a)$, as
\[w(a)=a_0+a_1+\cdots +a_{n-1}.\]
For integers $a$ divisible by $q-1$, we define $w(a)=0$.

Next let $\xi_{q-1}$ be a complex primitive $(q-1)$th root of unity.
Fix any prime ideal $\mathfrak{p}$ in $\Zz[\xi_{q-1}]$ lying over
$p$. Then $\Zz[\xi_{q-1}]/\mathfrak{p}$ is a finite field of order
$q$, which we identify with $\Ff_q$. Let $\omega_{\mathfrak{p}}$ be
the Teichm\"uller character on $\Ff_q$, i.e., an isomorphism
$$\omega_{\mathfrak{p}}: \Ff_q^*\rightarrow
\{1,\xi_{q-1},\xi_{q-1}^2,\dots ,\xi_{q-1}^{q-2}\}$$ satisfying
\begin{equation}\label{eq2.3}
\omega_{\mathfrak{p}}(\alpha)\quad ({\rm
mod}\hspace{0.1in}{\mathfrak{p}})=\alpha,
\end{equation}
for all $\alpha$ in $\Ff_q^*$. The Teichm\"uller character
$\omega_{\mathfrak{p}}$ has order $q-1$; hence it generates all
multiplicative characters of $\Ff_q$.

Let $\mathfrak{P}$ be the prime ideal of $\Zz[\xi_{q-1},\xi_p]$
lying above $\mathfrak{p}$. For an integer $a$, let $\nu
_{\mathfrak{P}}(g(\omega_{\mathfrak p}^{-a}))$ denote the
$\mathfrak{P}$-adic valuation of $g(\omega_{\mathfrak p}^{-a})$. The
following classical theorem is due to Stickelberger (see
\cite[p.~7]{lang}, \cite[p.~344]{be1}).

 \begin{teor}\label{stick}
Let $p$ be a prime, and $q=p^n$. Let $a$ be any integer not
divisible by $q-1$. Then
$$\nu_{\mathfrak{P}}(g(\omega_{\mathfrak p}^{-a}))=w(a).$$
\end{teor}

Next we state Hou's lemma using the notation developed in this
paper.

\begin{lem} {\rm (\cite{hou})} \label{hou1}
Let $f: \Ff_{3^n}\rightarrow \Ff_3$ be a function. We have\

{\rm (i)} $f$ is a ternary bent function if and only if
$\nu_3(S_f(b))=\frac{n}{2}$ for all $b\in\Ff_{3^n}$.\

{\rm (ii)} $f$ is a weakly regular bent function if and only if
$\nu_3(S_f(0))=\frac{n}{2}$ and $\nu_3(S_f(b)-S_f(0))>\frac{n}{2}$
for all $b\in \Ff_{3^n}^*$.
\end{lem}

\section{Proofs of the Main Results}\label{proof}

We will first prove Theorem~\ref{minor}. The proof is relatively
easy since most of the work has been done in \cite{fl}.

Let $F:\Ff_{p^n}\rightarrow\Ff_{p^n}$ be a function, and
$\omega=e^{\frac{2\pi i}{p}}$ be a primitive complex $p$th root of
unity. In \cite{fl}, the following notation was introduced:
$$S_F(a,b)=\sum_{x \in \Ff_q} \omega^{\Tr_n(aF(x)+bx)}$$
and $$K=\mathbb{Q}(\omega),\; W_K^{+}=\{\omega^i\,|\,0 \le i \le
p-1\}, \; W_K^{-}=\{-\omega^i\,|\,0 \le i \le p-1\}.$$ Note that
$W_K=W_K^{+}\cup W_K^-$ is the group of roots of unity in $K^*$. We
quote the following theorem from \cite{fl}.

\begin{teor}{\em(\cite{fl})}\label{cm}
Let $q$ be an odd prime power. Let $F$ be a planar function on
$\Ff_q$ with $F(0)=0$ and $F(-x)=F(x)$ for all $x\in \Ff_q$. Then we
have

i) \begin{eqnarray*}
   \sum_{a \in \Ff_q^*}S_F(a, 0)&=&0\\
   \sum_{a, b \in \Ff_q}S_F(a,b)&=&\sum_{a, b\in \Ff_q}S_F(a, b)^2=q^2
\end{eqnarray*}
ii) For all $a \in \Ff_q^*$ and $b \in \Ff_q$
$$S_F(a, b)=\varepsilon_{a, b}(\sqrt{p^*})^n,\quad \varepsilon_{a, b}\in W_K,$$
where $p^*=(-1)^{\frac{p-1}{2}}p$. Moreover, if $F$ is of
Dembowski-Ostrom type or $F$ is the Coulter-Matthews planar
function, then
$$\varepsilon_{a, 0}\in \{\pm 1\}\quad\mbox{and}\quad \varepsilon_{a, b}\cdot\varepsilon_{a, 0}\in W_K^+.$$
\end{teor}

We are ready to give the proof of Theorem~\ref{minor}.

\bpft{\ref{minor}}
Let $F:
\Ff_{3^n}\rightarrow \Ff_{3^n}$ be defined by $F(x)=x^{\frac {3^k
+1} {2}},\;\forall x\in \Ff_{3^n}$. For any nonzero $a\in\Ff_{3^n}$,
let $f:\Ff_{3^n}\rightarrow \Ff_{3}$ be defined by
$f(x)=\Tr_n(ax^{\frac {3^k +1} {2}}),\;\forall x\in \Ff_{3^n}$. By
Lemma~\ref{hou1}, it suffices to show that $\nu_3(S_f(0))=n/2$, and
for every $b\in \Ff_{3^n}^*$, $\nu_{3}(S_f(b)-S_f(0))>\frac{n}{2}$.

As $F$ is a planar function on $\Ff_{3^n}$, by Theorem~\ref{cm},
$$S_f(0)=S_F(a,0)=\varepsilon_{a,0}(\sqrt{-3})^n.$$
Therefore $\nu_3(S_f(0))=\frac{n}{2}$.

For any $b\in\Ff_{3^n}^*$, we have
$$S_f(b)-S_f(0)=\sum_{x \in \Ff_q}\omega^{\Tr_n(aF(x)-bx)}-\sum_{x \in
\Ff_q}\omega^{\Tr_n(aF(x))}=S_F(a,-b)-S_F(a,0).$$ By
Theorem~\ref{cm}, we have
$$S_F(a,-b)=\varepsilon_{a, -b}(\sqrt{-3})^n,\;
S_F(a,0)=\varepsilon_{a,0}(\sqrt{-3})^n,$$ and $$\varepsilon_{a,
0}\in \{\pm 1\}\quad\mbox{and}\quad \varepsilon_{a,
-b}\cdot\varepsilon_{a, 0}\in W_K^+.$$ Therefore,
\begin{eqnarray*}
S_f(b)-S_f(0)&=&(\sqrt{-3})^n(\varepsilon_{a, -b}-\varepsilon_{a,
0})\\
&=&(\sqrt{-3})^n\varepsilon_{a, 0}(\omega^j-1),
\end{eqnarray*}
where $\omega$ is a complex primitive cubic root of unity, and
$j\in\{0,1,2\}$.

Fix any prime ideal $\mathfrak{p}$ in $\Zz[\xi_{q-1}]$ lying over
$3$. Let $\mathfrak{P}$ be the prime ideal of
$\Zz[\xi_{q-1},\omega]$ lying above $\mathfrak{p}$. Since $\nu
_{\mathfrak{P}}(3)=2$, we see that
$$\nu_3(S_f(b)-S_f(0))>\frac{n}{2} \iff \nu_{\mathfrak{P}}(S_f(b)-S_f(0))>n.$$
Note that for $j=0$, we have $\nu_{\mathfrak P}(\omega^j-1)=\infty$;
and for $j=1$ or 2, we have  $\nu_{\mathfrak P}(\omega^j-1)=1$. As
$\nu_{\mathfrak P}(\sqrt{-3})^n=n$, we have
$$\nu_{\mathfrak P}(S_f(b)-S_f(0))=\nu_{\mathfrak P}(\omega^j-1)+\nu_{\mathfrak P}(\sqrt{-3})^n>n.$$
Hence we have shown that $\nu_3(S_f(b)-S_f(0))>\frac{n}{2}$. The
proof of theorem is now complete. \epft
\begin{rem}
It was shown in \cite{hou} that if $f:\Ff_p^n\rightarrow\Ff_p$ is a
weakly regular bent function and $(p-1)n\geq 4$, then
\begin{equation}\label{degbound}{\rm deg}(f)\leq \frac {(p-1)n}{2}.\end{equation}
In \cite{hou}, after the proof of the bound in (\ref{degbound}), it
was mentioned that when $p$ and $n$ are both odd with $n\geq 3$, it
is not known if the bound in (\ref{degbound}) is attainable. Let
$n\geq 3$ be an integer. Then the Coulter-Matthews bent functions
$\Tr_n(ax^{\frac {3^{n-1} +1} {2}})$ has degree $n$. Therefore by
Theorem~\ref{minor}, these functions provide examples of weakly
regular bent functions of degrees attaining the bound in
(\ref{degbound}).
\end{rem}

We now make some preparation for the proof of Theorem~\ref{major}.
Let $C_i$ $(i=0,1,2,3)$ denote the {\em cyclotomic classes} of order
four in the multiplicative group of $\Ff_{p^n}$, i.e.,
$C_i=\{\xi^{4t+i}\ |\ t=0,\dots,f-1\}$, where $\xi$ is a primitive
element of $\Ff_{p^n}$ and $f=\frac{p^n-1}{4}$. Throughout this
section all expressions in the indices numbering the cyclotomic
classes are taken modulo $4$.

\begin{lem} {\em(\cite{hk2})}\label{simplecase}
Let $p$ be an odd prime with $p\equiv 3\ (\bmod\;4)$ and let $n=2k$
with $k$ odd. Raising elements of $C_i$ to the $(p^k+1)^{\rm th}$
power results in a $\frac{p^k+1}{2}\;$-to-$1$ mapping onto the
cyclotomic classes of order two in the multiplicative group of
$\Ff_{p^k}$. Moreover, $C_0$ and $C_2$ map onto the squares and
$C_1$ and $C_3$ onto the non-squares in $\Ff_{p^k}^*$.
\end{lem}

\bpf Take the following polynomial over $\Ff_p$ that factors in
$\Ff_{p^k}$ as
\[p(z)=z^{\frac{p^n-1}{4}}-1=(z^t)^{p^k-1}-1=\prod_{\alpha\in\Ff_{p^k}^*}(z^t-\alpha)\]
where $t=\frac{p^k+1}{4}$. The roots of $p(z)$ are exactly all the
elements from $C_0$. Therefore, it can be concluded that raising
elements of $C_0$ to the power of $t$ results in a $t$-to-$1$
mapping onto the multiplicative group of $\Ff_{p^k}$. In general,
raising elements of $C_i=\xi^i C_0$ to the $t^{\rm th}$ power
results in a $t$-to-$1$ mapping onto the coset
$\xi^{it}\Ff_{p^k}^*$.

Let $\eta=\xi^{p^k+1}$ be a primitive element of $\Ff_{p^k}$. When
$k$ is odd, the cyclic subgroups generated by $\eta^2$ and by
$\eta^4$ are equal since they have the same multiplicative order
equal to
\[{\rm ord}\,(\eta^4)=\frac{p^k-1}{\gcd(p^k-1,4)}=\frac{p^k-1}{2}={\rm
ord}\,(\eta^2)\enspace.\] Thus, raising elements of $\Ff_{p^k}^*$ to
the fourth power is a mapping onto the subgroup generated by
$\eta^2$ and since both $\alpha$ and $-\alpha$ produce the same
image for any $\alpha\in\Ff_{p^k}^*$, this is a $2$-to-$1$ mapping.

Also note that $\xi^{4it}=\xi^{i(p^k+1)}=\eta^i$. Therefore,
combination of these two mappings that is equivalent to raising
elements of $C_i$ to the power of $4t=p^k+1$, results in a
$\frac{p^k+1}{2}\;$-to-$1$ mapping onto the cyclotomic classes of
order two in $\Ff_{p^k}^*$. Moreover, $C_0$ and $C_2$ map onto the
squares and $C_1$ and $C_3$ onto the non-squares in $\Ff_{p^k}^*$.
\epf

\begin{lem} {\em(\cite{hk2})}
 \label{le:1}
Let $p$ be an odd prime with $p\equiv 3\ (\bmod\;8)$ and let $n=2k$
with $k$ odd. Then for any $c\in\Ff_{p^k}^*$ and $z\in\Ff_{p^n}^*$, and any cyclotomic class $C_j$ $(j=0,1,2,3)$
\[\sum_{y\in C_j}\omega^{{\rm Tr}_n\left(cz^{p^k}y\right)}=
\left\{\begin{array}{ll}
\frac{3p^k-1}{4},&\ \mbox{if}\quad z\in C_{j+2}\\
-\frac{p^k+1}{4},&\ \mbox{otherwise}\ .
\end{array}\right.\]
\end{lem}

This lemma is a direct consequence of part (1) of the following
general theorem \cite{mye} on uniform cyclotomy. See also
\cite{bmw}.

\begin{teor}{\em (\cite{mye})}
Let $q=p^{n}$ be a prime power, let $e>1$ be a divisor of $q-1$ and
let $C_i$, $0\leq i\leq e-1$, be the cyclotomic classes of order
$e$. Assume there exists a positive integer $j$ such that $p^j\equiv
-1\; (mod$ $e)$, and assume $j$ is the smallest such integer. Moreover
assume that $n=2j\gamma$. Then the cyclotomic periods
$\eta_i=\sum_{x\in C_i}\omega^{{\rm Tr}_n(x)}$ are given as follows:

(1) If $\gamma, p, \frac {p^j+1} {e}$ are all odd, then
$$\eta_{\frac {e} {2}}=\frac {(e-1)p^{j\gamma}-1} {e},\; \eta_i=\frac {-1-p^{j\gamma}} {e}, i\neq e/2.$$

(2) In all other cases
$$\eta_0=\frac {-1-(-1)^{\gamma}(e-1)p^{j\gamma}} {e},\; \eta_i=\frac {(-1)^{\gamma}p^{j\gamma}-1} {e}, i\neq 0.$$
\end{teor}

\begin{lem}{\em (\cite{hk2})}\label{conjugate}
Let $p$ be an odd prime with $p\equiv 3\ (\bmod\;8)$ and let $n=2k$
with $k$ odd. For any $c\in\Ff_{p^k}$ and $j=0,1,2,3$ denote
\[T_j=\sum_{x\in C_j}\omega^{{\rm Tr}_k\left(c(x+1)^{p^k+1}-c\right)}\enspace.\]
Then for any $j$
\[-\overline{T_j}=\omega^{{\rm
Tr}_k(c)}T_{j+2}+\frac{p^k+1}{4}\big(\omega^{{\rm Tr}_k(c)}+1\big)\]
where the bar over a complex value denotes the complex conjugate and
the indices are taken modulo $4$.
\end{lem}

\bpf First, it is easy to see that for any nonzero $c\in\Ff_{p^k}$
\begin{eqnarray}
 \label{eq:12}
\nonumber&&1+T_0+T_1+T_2+T_3=\sum_{x\in\Ff_{p^n}}\omega^{{\rm
Tr}_k\left(c(x+1)^{p^k+1}-c\right)}\\
\nonumber&&=\omega^{-{\rm
Tr}_k(c)}\sum_{y\in\Ff_{p^n}}\omega^{{\rm
Tr}_k\left(cy^{p^k+1}\right)}\stackrel{(\ast)}{=}\omega^{-{\rm
Tr}_k(c)}\left((p^k+1)\sum_{z\in\Ff_{p^k}^*}\omega^{{\rm
Tr}_k(cz)}+1\right)\\
&&=-p^k\omega^{-{\rm Tr}_k(c)}
\end{eqnarray}
where $(\ast)$ holds since raising elements of $\Ff_{p^n}^*$ to the
$(p^k+1)^{\rm th}$ power is a $(p^k+1)$-to-$1$ mapping onto
$\Ff_{p^k}^*$ as proved in \cite[Lemma~1]{DeGo69}.

Let $C_i\cdot C_j$ denote the {\em strong union} of $C_i$ and $C_j$,
i.e., the set of elements of $\Ff_{p^n}$ that can be represented as
a sum of two addends from $C_i$ and $C_j$, respectively, and
counting the multiplicity of such a representation. Thus, $C_i\cdot
C_j$ consists of the elements
$\xi^{4t+i}+\xi^{4d+j}=\xi^{4d+j}(1+\xi^{4(t-d)+i-j})$ for all
$t,d=0,\dots,f-1$. Therefore,
\begin{eqnarray}
 \label{eq:15}
\nonumber C_i\cdot C_j&=&C_j(1+C_{i-j})\\
\nonumber&=&(i-j,0)C_j\cup(i-j,1)C_{j+1}\cup(i-j,2)C_{j+2}\cup(i-j,3)C_{j+3}\\
&=&(i-j,-j)C_0\cup(i-j,1-j)C_1\cup(i-j,2-j)C_2\cup(i-j,3-j)C_3\enspace\enspace\enspace
\end{eqnarray}
if $i\neq j$ and otherwise, since $-1\in C_0$,
\begin{equation}
 \label{eq:16}
C_i\cdot C_i=(0,-i)C_0\cup(0,1-i)C_1\cup(0,2-i)C_2\cup(0,3-i)C_3\cup
f\{0\}
\end{equation}
where $(i,j)$ denotes the {\em cyclotomic number} that is equal to
the number of elements $x\in C_i$ such that $x+1\in C_j$ and
$f\{0\}$ denotes the zero-element of $\Ff_{p^n}$ taken with the
multiplicity $f$. The components $i$, $j$ in cyclotomic numbers are
taken modulo $4$.

Also denote $C_i^j=\{x\in C_i\ |\ 1+x\in C_j\}$ (obviously,
$|C_i^j|=(i,j)$). In our case $-1\in C_0$ and we can prove that
$(i,j)=(j,i)$. Indeed, the elements of $C_i^j$ correspond to the
pairs $(t,d)$ with $t,d\in\{0,\dots,f-1\}$ that satisfy the equation
$\xi^{4t+i}+1=\xi^{4d+j}$. Multiplying both sides of the equation by
$-1=\xi^{4l}$ we get the equivalent equation
$\xi^{4(d+l)+j}+1=\xi^{4(t+l)+i}$ whose solutions give the elements
of $C_j^i$. Therefore, for any $i\in\{0,1,2,3\}$ we have
\[\sum_{j=0}^3(j,i)=\sum_{j=0}^3(i,j)=|C_i^0\cup C_i^1\cup C_i^2\cup
C_i^3|=\left\{\begin{array}{ll}
|C_i|=f,&\ \mbox{if}\quad i\neq 0\\
|C_0\backslash\{-1\}|=f-1,&\ \mbox{otherwise}
\end{array}\right.\]
since $-1+1=0$ that does not belong to any $C_i$. A good
introduction into this subject can be found in \cite{St67}. Now for
$i,j\in\{0,1,2,3\}$ and $i\neq j$ we evaluate the product
\begin{eqnarray}
 \label{eq:13}
\nonumber T_i\overline{T_j}&=&\sum_{x\in C_i,\,y\in C_j}\omega^{{\rm
Tr}_k\left(c(x+1)^{p^k+1}-c-c(y+1)^{p^k+1}+c\right)}\\
\nonumber&=&\sum_{x\in C_i,\,y\in C_j}\omega^{{\rm
Tr}_k\left(c\left(x^{p^k+1}-y^{p^k+1}+(x-y)^{p^k}+(x-y)\right)\right)}\\
\nonumber&\stackrel{\{C_j=-C_j\}}{=}&\sum_{x\in C_i,\,y\in
C_j}\omega^{{\rm
Tr}_k\left(c\left(x^{p^k+1}-y^{p^k+1}+(x+y)^{p^k}+(x+y)\right)\right)}\\
\nonumber&=&\sum_{z\in C_i\cdot C_j}\omega^{{\rm
Tr}_k\left(c\left((z-y)^{p^k+1}-y^{p^k+1}+z^{p^k}+z\right)\right)}\\
\nonumber&=&\sum_{z\in C_i\cdot C_j}\omega^{{\rm
Tr}_k\left(c(z+1)^{p^k+1}-c\right)}\omega^{-{\rm
Tr}_k\left(c\left(zy^{p^k}+z^{p^k}y\right)\right)}\\
\nonumber&=&\sum_{z\in C_i\cdot
C_j}\omega^{{\rm Tr}_k\left(c(z+1)^{p^k+1}-c\right)}\omega^{-{\rm Tr}_n\left(cz^{p^k}y\right)}\\
&=&\sum_{t=0}^3\sum_{z\in C_t}\omega^{{\rm
Tr}_k\left(c(z+1)^{p^k+1}-c\right)}\sum_{r\in
C_{i-j}^{t-j}}\omega^{-{\rm Tr}_n\left(cz^{p^k}\frac{z}{1+r}\right)}
\end{eqnarray}
where $z=x+y\in C_i\cdot C_j$ and the value of $y$ is uniquely
defined by $z$. Therefore, if $z=x+y\in C_t$ with $x\in C_i$ and
$y\in C_j$ then $z=y(1+xy^{-1})$ with $xy^{-1}\in C_{i-j}^{t-j}$. By
(\ref{eq:15}), the multiplicity of $z\in C_t$ in $C_i\cdot C_j$ is
equal to $(i-j,t-j)=\big|C_{i-j}^{t-j}\big|$. Thus, for a fixed
$z\in C_t$ the set $\big\{\frac{z}{1+r}\ |\ r\in
C_{i-j}^{t-j}\big\}$ contains all $(i-j,t-j)$ values for $y\in C_j$
that correspond to this $z$ taken with the appropriate multiplicity
$(i-j,t-j)$ as a member of $C_i\cdot C_j$. For $i=j$ we just have
additionally to consider the zero-element of $\Ff_{p^n}$ that is
found in $C_i\cdot C_i$ with the multiplicity $f$ (see
(\ref{eq:16})). Then
\begin{equation}
 \label{eq:17}
T_i\overline{T_i}=\sum_{t=0}^3\sum_{z\in C_t}\omega^{{\rm
Tr}_k\left(c(z+1)^{p^k+1}-c\right)}\sum_{r\in
C_0^{t-i}}\omega^{-{\rm
Tr}_n\left(cz^{p^k}\frac{z}{1+r}\right)}+f\enspace.
\end{equation}

Let $t,j\in\{0,1,2,3\}$ and $z\in C_t$ be fixed. Then for any
$i\in\{0,1,2,3\}$ and $r\in C_{i-j}^{t-j}$ we have $\frac{z}{1+r}\in
C_j$. Further,
$\sum_{i=0}^3\big|C_{i-j}^{t-j}\big|=\sum_{i=0}^3(i,t-j)$ is equal
to $|C_{t-j}|=f$ if $t\neq j$ and is equal to $|C_0|-1=f-1$
otherwise. Since the cardinality of $C_j$ is $f$, we have proven
that
\[\left\{\frac{z}{1+r}\ |\ r\in C_{i-j}^{t-j},\
i=0,1,2,3\right\}=\left\{\begin{array}{ll}
C_j,&\ \mbox{if}\quad t\neq j\\
C_j\backslash\{z\},&\ \mbox{otherwise}
\end{array}\right.\]
since $r\neq 0$. Therefore, for any $t,j\in\{0,1,2,3\}$ and $z\in
C_t$
\begin{equation}
 \label{eq:24}
\sum_{i=0}^3\sum_{r\in C_{i-j}^{t-j}}\omega^{-{\rm
Tr}_n\left(cz^{p^k}\frac{z}{1+r}\right)}=\left\{\begin{array}{ll}
\sum_{y\in C_j}\omega^{-{\rm Tr}_n\left(cz^{p^k}y\right)},&\ \mbox{if}\quad t\neq j\\
\sum_{y\in C_j\backslash\{z\}}\omega^{-{\rm
Tr}_n\left(cz^{p^k}y\right)},&\ \mbox{otherwise}\,.
\end{array}\right.
\end{equation}

Note that since $n=2k$ and $p\equiv 3\ (\bmod\;8)$ then
$(p^n-1)/2\equiv 0\ (\bmod\;4)$ and $-1=\xi^{\frac{p^n-1}{2}}\in
C_0$. Therefore, $-C_j=C_j$ and
\begin{equation}
 \label{eq:25}
\overline{T_j}=\sum_{z\in C_j}\omega^{{\rm
Tr}_k\left(c(-z^{p^k+1}-z^{p^k}-z)\right)}=\sum_{z\in
C_j}\omega^{{\rm Tr}_k\left(c(-z^{p^k+1}+z^{p^k}+z)\right)}\enspace.
\end{equation}

Making use of Lemma~\ref{le:1} we get that
\begin{eqnarray*}
&&(T_0+T_1+T_2+T_3)\overline{T_j}\\
&&\stackrel{(\ref{eq:13},\ref{eq:17})}{=}\sum_{t=0}^3\sum_{z\in
C_t}\omega^{{\rm
Tr}_k\left(c(z+1)^{p^k+1}-c\right)}\sum_{i=0}^3\sum_{r\in
C_{i-j}^{t-j}}\omega^{-{\rm
Tr}_n\left(cz^{p^k}\frac{z}{1+r}\right)}+f\\
&&\stackrel{(\ref{eq:24})}{=}\sum_{t=0}^3\sum_{z\in C_t}\omega^{{\rm
Tr}_k\left(c(z+1)^{p^k+1}-c\right)}\sum_{y\in C_j}\omega^{-{\rm
Tr}_n\left(cz^{p^k}y\right)}\\
&&-\sum_{z\in C_j}\omega^{{\rm
Tr}_k\left(c(z+1)^{p^k+1}-c\right)}\omega^{-{\rm
Tr}_n\left(cz^{p^k+1}\right)}+f\\
&&=-\frac{p^k+1}{4}\sum_{t\neq
j+2}T_t+\frac{3p^k-1}{4}T_{j+2}-\sum_{z\in C_j}\omega^{{\rm
Tr}_k\left(c(-z^{p^k+1}+z^{p^k}+z)\right)}+f\\
&&\stackrel{(\ref{eq:25})}{=}-\frac{p^k+1}{4}(T_0+T_1+T_2+T_3)+p^k
T_{j+2}-\overline{T_j}+f\enspace.
\end{eqnarray*}
Now, using (\ref{eq:12}), we get
\begin{eqnarray*}
-p^k\omega^{-{\rm Tr}_k(c)}\overline{T_j}&=&(p^k\omega^{-{\rm
Tr}_k(c)}+1)\frac{p^k+1}{4}+p^k
T_{j+2}+\frac{p^n-1}{4}\quad\mbox{and}\\
-\overline{T_j}&=&\omega^{{\rm
Tr}_k(c)}T_{j+2}+\frac{p^k+1}{4}(\omega^{{\rm Tr}_k(c)}+1)
\end{eqnarray*}
that was claimed.
\epf

\begin{lem}{\em (\cite{hk1,HeKh07})}\label{kasami}
Let $n=2k$ and $a\in\Ff_{p^n}$ for an odd prime $p$. Then the
function $f$ defined by $f(x)={\rm Tr}_n\big(a x^{p^k+1}\big)\;
\forall x\in\Ff_{p^n}$ is bent if and only if $a+a^{p^k}\neq 0$.
Moreover, if the latter condition holds then $f$ is weakly regular
and for $b\in\Ff_{p^n}$, the corresponding Walsh transform
coefficient of $f$ is equal to
\[S_a(b)=-p^k\omega^{-{\rm
Tr}_k\left(\frac{b^{p^k+1}}{a+a^{p^k}}\right)}\enspace.\]
\end{lem}

For the proof of this lemma, we refer the reader to
(\cite{hk1,HeKh07}).

\bpft{\ref{major}} By Lemma \ref{hou1}, it suffices to show that
$\nu_3(S_a(0))=k$, and for every $b \in\Ff_{3^n}^*$,
$\nu_3(S_a(b)-S_a(0))>\frac{n}{2}$. First we will compute $S_a(0)$
and $S_a(b)-S_a(0)$.

Let $I=\xi^{\frac{3^n-1}{4}}$, where $I$ is a primitive $4^{\rm th}$
root of unity in $\Ff_{3^n}$ (obviously $I^2=-1$). As before, let
$C_i$, $0\leq i\leq 3$, be the cyclotomic classes of order 4 of
$\Ff_{3^n}$. Then any $x \in C_i$ satisfies
$x^{\frac{3^n-1}{4}}=\xi^{\frac{i(3^n-1)}{4}}=I^i$. Also
$a^{3^k}=aI$ and $\Tr_k^n(a)=a+a^{3^k}=a(I+1)$. On the other hand,
$\Tr_k^n(aI)=aI-a^{3^k}I=aI+a=a(I+1)=\Tr_k^n(a)$ since $3^k\equiv 3
\pmod 4$ for odd $k$. Therefore,

\begin{eqnarray}\label{sabCOMP}
S_a(b)-1&=&\sum_{x \in
\Ff_{3^n}}\omega^{\Tr_n(ax^{\frac{3^n-1}{4}+3^k+1}-bx)}-1
=\sum_{i=0}^{3}\sum_{x \in C_i}\omega^{\Tr_n(aI^i x^{3^k+1}-bx)}\notag\\
&=&\sum_{x \in C_0\cup C_1}\omega^{\Tr_k(a_1 x^{3^k+1}-bx-b^{3^k}x^{3^k})}+
\sum_{x \in C_2\cup C_3}\omega^{\Tr_k(-a_1 x^{3^k+1}-bx-b^{3^k}x^{3^k})}\notag\\
&=&\sum_{x \in C_0 \cup
C_1}\omega^{\Tr_k(a_1(x-\beta)^{3^k+1}-a_1\beta^{3^k+1})} +\sum_{x
\in C_2\cup
C_3}\omega^{-\Tr_k(a_1(x+\beta)^{3^k+1}-a_1\beta^{3^k+1})},
\end{eqnarray}
where $a_1=a(I+1)\neq 0$ belongs to $\Ff_{3^k}$ and $b=a_1
\beta^{3^k}$.

If $b=0$, then $\beta=0$. Using Lemma \ref{simplecase}, we have
\begin{eqnarray*}
S_a(0)&=&1+\sum_{x \in C_0\cup C_1}\omega^{\Tr_k(a_1 x^{3^k+1})}
+\sum_{x \in C_2\cup C_3}\omega^{-\Tr_k(a_1 x^{3^k+1})}\\
&=&1+\frac{3^k+1}{2}\sum_{y \in\Ff_{3^k}^*}(\omega^{\Tr_k(a_1
y)}+\omega^{-\Tr_k(a_1 y)})=-3^k.
\end{eqnarray*}
Therefore $\nu_3(S_a(0))=k=n/2$.

Next suppose $b\neq 0$. Then $\beta \neq 0$. Let
$c=a_1\beta^{3^k+1}$. We have $c\in\Ff_{3^k}^*$. Assuming that
$\beta^{-1}\in C_j$ (i.e., ind($\beta^{-1}$)$\equiv j\pmod 4$), we
have $\beta^{-1}C_i=C_{i+j}$ for any $i \in \{0, 1, 2, 3\}$. Now
making the substitution $x=\beta y$ in~(\ref{sabCOMP}), we have

\begin{eqnarray*}
S_a(b)-1&=&\sum_{y\in C_j\cup
C_{j+1}}\omega^{\Tr_k(c(y-1)^{3^k+1}-c)} +\sum_{y \in C_{j+2}\cup
C_{j+3}}\omega^{-\Tr_{k}(c(y+1)^{3^k+1}-c)}.
\end{eqnarray*}

Since $n=2k$, $(3^n-1)/2\equiv 0 \pmod 4$ and
$-1=\xi^{\frac{3^n-1}{2}}\in C_0$. Therefore, $-C_i=C_i$ and
$$\sum_{y \in C_i}\omega^{\Tr_k(c(y-1)^{3^k+1})}=\sum_{y \in C_i}\omega^{\Tr_k(c(y+1)^{3^k+1})}.$$
Let $T_i$ $(i=0,1,2,3)$ be defined as in Lemma \ref{conjugate}. Then
we have
\begin{equation*}
S_a(b)=1+T_j+T_{j+1}+\overline{T_{j+2}}+\overline{T_{j+3}},
\end{equation*}
where the bars denote complex conjugation. By Lemma \ref{conjugate}
we have
\begin{equation}\label{sum}
S_a(b)=1+T_j+T_{j+1}+\overline{T_{j+2}}+\overline{T_{j+3}}=\left(1-\omega^{\Tr_k(c)}\right)\left(T_j+T_{j+1}+\frac{3^k+1}{2}\right)-3^k,
\end{equation}
where
$$T_j=\sum_{x \in C_j}\omega^{\Tr_k(c(x+1)^{3^k+1}-c)}=\sum_{x \in C_j}\omega^{\Tr_n(2cx^{3^k+1}+cx)}
=\sum_{x \in C_j}\omega^{\Tr_n(-cx^{3^k+1}+cx)}, \;j=0,1,2,3.$$ Let
$\eta$ be a multiplicative character of $\Ff_{3^n}$ of order $4$.
Then
\begin{eqnarray*}
T_0&=&\sum_{x \in \Ff_{3^n}^*}\omega^{\Tr_n(-cx^{3^k+1}+cx)}\frac{1+\eta(x)+\eta^2(x)+\eta^3(x)}{4};\\
T_1&=&\sum_{x \in \Ff_{3^n}^*}\omega^{\Tr_n(-cx^{3^k+1}+cx)}\frac{1-i\eta(x)-\eta^2(x)+i\eta^3(x)}{4};\\
T_2&=&\sum_{x \in \Ff_{3^n}^*}\omega^{\Tr_n(-cx^{3^k+1}+cx)}\frac{1-\eta(x)+\eta^2(x)-\eta^3(x)}{4};\\
T_3&=&\sum_{x \in
\Ff_{3^n}^*}\omega^{\Tr_n(-cx^{3^k+1}+cx)}\frac{1+i\eta(x)-\eta^2(x)-i\eta^3(x)}{4}.
\end{eqnarray*}
So
\begin{eqnarray}
T_0+T_1&=&\frac{1}{2}\sum_{x \in \Ff_{3^n}^*}\omega^{\Tr_n(-cx^{3^k+1}+cx)}\nonumber\\
&+&\frac{1-i}{4}\sum_{x \in \Ff_{3^n}^*}\omega^{\Tr_n(-cx^{3^k+1}+cx)}\eta(x)\nonumber\\
&+&\frac{1+i}{4}\sum_{x \in
\Ff_{3^n}^*}\omega^{\Tr_n(-cx^{3^k+1}+cx)}\eta^3(x); \label{tsum0}
\end{eqnarray}

\begin{eqnarray}
T_1+T_2&=&\frac{1}{2}\sum_{x \in \Ff_{3^n}^*}\omega^{\Tr_n(-cx^{3^k+1}+cx)}\nonumber \\
&+&\frac{-i-1}{4}\sum_{x \in \Ff_{3^n}^*}\omega^{\Tr_n(-cx^{3^k+1}+cx)}\eta(x)\nonumber \\
&+&\frac{i-1}{4}\sum_{x \in
\Ff_{3^n}^*}\omega^{\Tr_n(-cx^{3^k+1}+cx)}\eta^3(x);\label{tsum1}
\end{eqnarray}

\begin{eqnarray}
T_2+T_3&=&\frac{1}{2}\sum_{x \in \Ff_{3^n}^*}\omega^{\Tr_n(-cx^{3^k+1}+cx)}\nonumber \\
&+&\frac{-1+i}{4}\sum_{x \in \Ff_{3^n}^*}\omega^{\Tr_n(-cx^{3^k+1}+cx)}\eta(x)\nonumber \\
&+&\frac{-i-1}{4}\sum_{x \in
\Ff_{3^n}^*}\omega^{\Tr_n(-cx^{3^k+1}+cx)}\eta^3(x);\label{tsum2}
\end{eqnarray}

\begin{eqnarray}
T_3+T_0&=&\frac{1}{2}\sum_{x \in \Ff_{3^n}^*}\omega^{\Tr_n(-cx^{3^k+1}+cx)}\nonumber \\
&+&\frac{1+i}{4}\sum_{x \in \Ff_{3^n}^*}\omega^{\Tr_n(-cx^{3^k+1}+cx)}\eta(x)\nonumber \\
&+&\frac{1-i}{4}\sum_{x \in
\Ff_{3^n}^*}\omega^{\Tr_n(-cx^{3^k+1}+cx)}\eta^3(x).\label{tsum3}
\end{eqnarray}

In the following, we will compute
$\sum_{x\in\Ff_{3^n}^*}\omega^{\Tr_n(-cx^{3^k+1}+cx)}$,
$\sum_{x\in \Ff_{3^n}^*}\omega^{\Tr_n(-cx^{3^k+1}+cx)}\eta(x)$,
and
$\sum_{x \in\Ff_{3^n}^*}\omega^{\Tr_n(-cx^{3^k+1}+cx)}\eta^3(x)$, respectively.

By Lemma \ref{kasami}, we have
\begin{eqnarray}
&&\sum_{x \in \Ff_{3^n}^*}\omega^{\Tr_n(-cx^{3^k+1}+cx)}=\sum_{x \in \Ff_{3^n}}\omega^{\Tr_n(-cx^{3^k+1}+cx)}-\omega^{\Tr_n(0)} \nonumber \\
&&=-3^k\omega^{-\Tr_k(c)}-1=-3^k\omega^{\Tr_n(c)}-1. \label{firsum}
\end{eqnarray}
Next we will compute $\sum_{x
\in\Ff_{3^n}^*}\omega^{\Tr_n(-cx^{3^k+1}+cx)}\eta(x)$.\

To simplify notation we write $L=\Ff_{3^n}$, and
$g_{L}(\chi)=\sum_{x \in \Ff_{3^n}^*}\chi(x)\omega^{\Tr_n(x)}$. Then
$$\omega^{\Tr_n(x)}=\frac{1}{3^n-1}\sum_{\chi \in \widehat{L^*}}g_{L}(\chi)\overline{\chi}(x).$$
With this notation, we have
\begin{eqnarray*}
&&\sum_{x \in L^*}\omega^{\Tr_n(-cx^{3^k+1}+cx)}\eta(x)=\sum_{x \in L^*}\eta(x)\omega^{\Tr_n(-cx^{3^k+1})}\omega^{\Tr_n(cx)}\\
&&=\sum_{x \in L^*}\eta(x)\frac{1}{3^n-1}\sum_{\chi_1 \in \widehat{L^*}}g_{L}(\chi_1)\overline{\chi_1}(-cx^{3^k+1})\frac{1}{3^n-1}\sum_{\chi_2 \in \widehat{L^*}}g_L(\chi_2)\overline{\chi_2}(cx)\\
&&=\frac{1}{(3^n-1)^2}\sum_{\chi_1}\sum_{\chi_2}g_L(\chi_1)g_L(\chi_2)\overline{\chi_1}(-c)\overline{\chi_2}(c)
\sum_{x \in L^*}\eta(x)\overline{\chi_1}(x^{3^k+1})\overline{\chi_2}(x)\\
&&=\frac{1}{(3^n-1)^2}\sum_{\chi_1}\sum_{\chi_2}g_L(\chi_1)g_L(\chi_2)\overline{\chi_1}(-c)\overline{\chi_2}(c)
\sum_{x \in
L^*}\overline{\chi_1}^{3^k+1}(x)\overline{\chi_2}(x)\eta(x).
\end{eqnarray*}
If $\chi_2=\overline{\chi_1}^{3^k+1}\eta$, then for any $x \in L^*$,
$\overline{\chi_1}^{3^k+1}(x)\overline{\chi_2}(x)\eta(x)=1$.
Otherwise
$$\sum_{x \in L^*}\overline{\chi_1}^{3^k+1}(x)\overline{\chi_2}(x)\eta(x)=0.$$
Thus,
\begin{eqnarray*}
\sum_{x \in L^*}\omega^{\Tr_n(-cx^{3^k+1}+cx)}\eta(x)&=&\frac{1}{3^n-1}\sum_{\chi_1}g_L(\chi_1)g_L(\overline{\chi_1}^{3^k+1}\eta)\overline{\chi_1}(-c)\chi_1^{3^k+1}(c)\overline{\eta}(c)\\
&=&\frac{\overline{\eta}(c)}{3^n-1}\sum_{\chi_1}g_L(\chi_1)g_L(\overline{\chi_1}^{3^k+1}\eta)\chi_1(-c).
\end{eqnarray*}
So
\begin{eqnarray}
&&\sum_{x \in L^*}\omega^{\Tr_n(-cx^{3^k+1}+cx)}\eta(x) \nonumber \\
&&=\frac{\overline{\eta}(c)}{3^n-1}\sum_{b=0}^{3^n-2}g_L(\omega_{\mathfrak{p}}^{-b})g_L(\omega_{\mathfrak{p}}^{(3^k+1)b+\frac{3^{2k}-1}{4}})\omega_{\mathfrak{p}}^{-b}(-c),\label{secsum}
\end{eqnarray}
where ${\mathfrak{p}}$ is a prime ideal in $\Zz[\xi_{q-1}]$ lying
above 3 and $\omega_{\mathfrak{p}}$ is the Teichm\"uller character
of $L$.

Similarly, we can compute $\sum_{x \in
L^*}\omega^{\Tr_n(-cx^{3^k+1}+cx)}\eta^3(x)$ as follows:
\begin{eqnarray}
&&\sum_{x \in L^*}\omega^{\Tr_n(-cx^{3^k+1}+cx)}\eta^3(x) \nonumber \\
&&=\frac{1}{(3^n-1)^2}\sum_{\chi_1}\sum_{\chi_2}g_L(\chi_1)g_L(\chi_2)\overline{\chi_1}(-c)\overline{\chi_2}(c)
\sum_{x \in L^*}\overline{\chi_1}^{3^k+1}(x)\overline{\chi_2}(x)\eta^3(x) \nonumber\\
&&=\frac{\overline{\eta^3}(c)}{3^n-1}\sum_{\chi_1}g_L(\chi_1)g_L(\overline{\chi_1}^{3^k+1}\eta^3)\chi_1(-c) \nonumber \\
&&=\frac{\overline{\eta^3}(c)}{3^n-1}\sum_{b=0}^{3^n-2}g_L(\omega_{\mathfrak{p}}^{-b})g_L(\omega_{\mathfrak{p}}^{(3^k+1)b+\frac{3(3^{2k}-1)}{4}})\omega_{\mathfrak{p}}^{-b}(-c).
\label{thirdsum}
\end{eqnarray}
If $\beta^{-1} \in C_0$, then by (\ref{sum}), (\ref{tsum0}),
(\ref{firsum}), (\ref{secsum}) and (\ref{thirdsum}), we have
\begin{eqnarray}
S_a(b)&=&(1-\omega^{\Tr_k(c)})(T_0+T_1+\frac{3^k+1}{2})-3^k \nonumber\\
&=&(1-\omega^{\Tr_k(c)})[-\frac{1}{2}3^k\omega^{\Tr_n(c)}+\frac{3^k}{2} \nonumber\\
&+&\frac{1-i}{4}\frac{\overline{\eta}(c)}{3^n-1}\sum_{b=0}^{3^n-2}g_L(\omega_{\mathfrak{p}}^{-b})g_L(\omega_{\mathfrak{p}}^{(3^k+1)b+\frac{3^{2k}-1}{4}})\omega_{\mathfrak{p}}^{-b}(-c) \nonumber\\
&+&\frac{1+i}{4}\frac{\overline{\eta^3}(c)}{3^n-1}\sum_{b=0}^{3^n-2}g_L(\omega_{\mathfrak{p}}^{-b})g_L(\omega_{\mathfrak{p}}^{(3^k+1)b+\frac{3(3^{2k}-1)}{4}})\omega_{\mathfrak{p}}^{-b}(-c)]-3^k.
\end{eqnarray} \label{walshform}

Since $S_a(0)=-3^k$, we have
\begin{eqnarray}
S_a(b)-S_a(0)&=&(1-\omega^{\Tr_k(c)})[-\frac{1}{2}3^k\omega^{\Tr_n(c)}+\frac{3^k}{2}\nonumber \\
&+&\frac{1-i}{4}\frac{\overline{\eta}(c)}{3^n-1}\sum_{b=0}^{3^n-2}g_L(\omega_{\mathfrak{p}}^{-b})g_L(\omega_{\mathfrak{p}}^{(3^k+1)b+\frac{3^{2k}-1}{4}})\omega_{\mathfrak{p}}^{-b}(-c)\nonumber \\
&+&\frac{1+i}{4}\frac{\overline{\eta^3}(c)}{3^n-1}\sum_{b=0}^{3^n-2}g_L(\omega_{\mathfrak{p}}^{-b})g_L(\omega_{\mathfrak{p}}^{(3^k+1)b+\frac{3(3^{2k}-1)}{4}})\omega_{\mathfrak{p}}^{-b}(-c)].\label{diff0}
\end{eqnarray}

Similarly, when $\beta^{-1} \in C_1$, we have
\begin{eqnarray}
S_a(b)-S_a(0)&=&(1-\omega^{\Tr_k(c)})[-\frac{1}{2}3^k\omega^{\Tr_n(c)}+\frac{3^k}{2}\nonumber \\
&+&\frac{-i-1}{4}\frac{\overline{\eta}(c)}{3^n-1}\sum_{b=0}^{3^n-2}g_L(\omega_{\mathfrak{p}}^{-b})g_L(\omega_{\mathfrak{p}}^{(3^k+1)b+\frac{3^{2k}-1}{4}})\omega_{\mathfrak{p}}^{-b}(-c)\nonumber \\
&+&\frac{i-1}{4}\frac{\overline{\eta^3}(c)}{3^n-1}\sum_{b=0}^{3^n-2}g_L(\omega_{\mathfrak{p}}^{-b})g_L(\omega_{\mathfrak{p}}^{(3^k+1)b+\frac{3(3^{2k}-1)}{4}})\omega_{\mathfrak{p}}^{-b}(-c)].\label{diff1}
\end{eqnarray}

when $\beta^{-1} \in C_2$, we have
\begin{eqnarray}
S_a(b)-S_a(0)&=&(1-\omega^{\Tr_k(c)})[-\frac{1}{2}3^k\omega^{\Tr_n(c)}+\frac{3^k}{2}\nonumber \\
&+&\frac{-1+i}{4}\frac{\overline{\eta}(c)}{3^n-1}\sum_{b=0}^{3^n-2}g_L(\omega_{\mathfrak{p}}^{-b})g_L(\omega_{\mathfrak{p}}^{(3^k+1)b+\frac{3^{2k}-1}{4}})\omega_{\mathfrak{p}}^{-b}(-c)\nonumber \\
&+&\frac{-i-1}{4}\frac{\overline{\eta^3}(c)}{3^n-1}\sum_{b=0}^{3^n-2}g_L(\omega_{\mathfrak{p}}^{-b})g_L(\omega_{\mathfrak{p}}^{(3^k+1)b+\frac{3(3^{2k}-1)}{4}})\omega_{\mathfrak{p}}^{-b}(-c)],\label{diff2}
\end{eqnarray}

and when $\beta^{-1} \in C_3$, we have
\begin{eqnarray}
S_a(b)-S_a(0)&=&(1-\omega^{\Tr_k(c)})[-\frac{1}{2}3^k\omega^{\Tr_n(c)}+\frac{3^k}{2}\nonumber \\
&+&\frac{1+i}{4}\frac{\overline{\eta}(c)}{3^n-1}\sum_{b=0}^{3^n-2}g_L(\omega_{\mathfrak{p}}^{-b})g_L(\omega_{\mathfrak{p}}^{(3^k+1)b+\frac{3^{2k}-1}{4}})\omega_{\mathfrak{p}}^{-b}(-c)\nonumber \\
&+&\frac{1-i}{4}\frac{\overline{\eta^3}(c)}{3^n-1}\sum_{b=0}^{3^n-2}g_L(\omega_{\mathfrak{p}}^{-b})g_L(\omega_{\mathfrak{p}}^{(3^k+1)b+\frac{3(3^{2k}-1)}{4}})\omega_{\mathfrak{p}}^{-b}(-c)].\label{diff3}
\end{eqnarray}

Let $\mathfrak{P}$ be the prime ideal of $\Zz[\xi_{q-1},\xi_3]$
lying above $\mathfrak{p}$. Since $\nu _{\mathfrak{P}}(3)=2$, we see
that
$$ \nu_3(S_a(b)-S_a(0))>\frac{n}{2} \iff \nu_{\mathfrak{P}}(S_a(b)-S_a(0))>n=2k.$$
Note that $\omega^{\Tr_k(c)}=1$, $\omega$ or $\omega^2$. Hence
$\nu_{\mathfrak{P}}(1-\omega^{\Tr_k(c)})=\infty$ or $1$. Using the
expressions of $S_a(b)-S_a(0)$ in (\ref{diff0}), (\ref{diff1}),
(\ref{diff2}), and (\ref{diff3}), we see that
$\nu_{\mathfrak{P}}(S_a(b)-S_a(0))>n$ if

\begin{equation}\label{finalequation1}
\nu_{\mathfrak{P}}\left(\sum_{b=0}^{3^n-2}g_L(\omega_{\mathfrak{p}}^{-b})g_L(\omega_{\mathfrak{p}}^{(3^k+1)b+\frac{3^{2k}-1}{4}})\omega_{\mathfrak{p}}^{-b}(-c)\right)
\ge 2k
\end{equation}
and
\begin{equation}\label{finalequation2}
\nu_{\mathfrak{P}}\left(\sum_{b=0}^{3^n-2}g_L(\omega_{\mathfrak{p}}^{-b})g_L(\omega_{\mathfrak{p}}^{(3^k+1)b+\frac{3(3^{2k}-1)}{4}})\omega_{\mathfrak{p}}^{-b}(-c)\right)\ge
2k.
\end{equation}
By Theorem~\ref{stick} and the fact that $g_L(\chi_0)=-1$, where
$\chi_0$ is the trivial multiplicative character of $\Ff_{3^n}$, we
have for any $b$, $0\leq b\leq 3^n-2$,
$$\nu_{\mathfrak{P}}\left(g_L(\omega_{\mathfrak{p}}^{-b})g_L(\omega_{\mathfrak{p}}^{(3^k+1)b+\frac{3^{2k}-1}{4}})\omega_{\mathfrak{p}}^{-b}(-c)\right)=w(b)+w\left(-(3^k+1)b-\frac{3^{2k}-1}{4}\right)$$
and
$$\nu_{\mathfrak{P}}\left(g_L(\omega_{\mathfrak{p}}^{-b})g_L(\omega_{\mathfrak{p}}^{(3^k+1)b+\frac{3(3^{2k}-1)}{4}})\omega_{\mathfrak{p}}^{-b}(-c)\right)=w(b)+w\left(-(3^k+1)b-\frac{3\left(3^{2k}-1\right)}{4}\right).$$
 Therefore if we
can prove that for each $b$, $0 \le b \le q-2$,
\begin{equation}\label{wtinequ1}
w(b)+w\left(-(3^k+1)b-\frac{3^{2k}-1}{4}\right)\ge 2k
\end{equation}
and
\begin{equation}\label{wtinequ2}
w(b)+w\left(-(3^k+1)b-\frac{3\left(3^{2k}-1\right)}{4}\right) \ge
2k,
\end{equation}
then $\nu_3(S_a(b)-S_a(0))>n/2$; and it follows that $f$ is weakly
regular bent by Lemma \ref{hou1}. We will give proofs of
(\ref{wtinequ1}) and (\ref{wtinequ2}) in the next two sections,
which will complete the proof of Theorem~\ref{major} \epft

\section{The $p$-ary modular add-with-carry algorithm}
In a sequence of papers \cite{gauss} \cite{hx} \cite{cyc} \cite{gmw}
a systematic method has been developed to derive binary weight
inequalities. Here we generalize this approach to $p$-ary weight
inequalities. As in the binary case, the idea is to analyze the
digit-wise contributions to the weights in the inequality using the
carries generated by a modular add-with-carry algorithm for the
$p$-ary numbers involved. Essentially, this approach enables the
analysis of the {\em global\/} properties of the weights in terms of
{\em local\/}, digit-wise contributions.

Then, these local contributions can be analyzed {\em for all word
lengths simultaneously\/} in a {\em finite\/} weighted directed
graph that models these local contributions. This graph has the
property that valid computations are in one-to-one correspondence
with directed closed walks in the graph. Hence the original weight
inequality gets transformed into a bound on the sum of the
arc-weights of directed closed walks in this graph as a function of
the length of the walk. In principle, such a bound can then be
verified by inspection, either directly (if the graph is
sufficiently small) or with the aid of a computer. Alternatively, a
detailed analysis of the properties of the graph, possibly with the
aid of a computer, can be used to devise a mathematical proof
(although such proofs can be quite tedious, see e.g.\ \cite{niho}).

We start with the derivation of the $p$-ary modular add-with-carry algorithm.
Our aim is to prove the following theorem.
\begin{teor}[Modular $p$-ary add-with-carry algorithm] \label{TAWCmod}
Let $a^{(1)}, \ldots, a^{(m)}$ be $m$ integers, and let
the integer $s$ satisfy
\[s \equiv t_1a^{(1)}+t_2a^{(2)}+\cdots +t_ma^{(m)} \bmod p^n-1\]
for nonzero integers $t_1, t_2, \ldots, t_m$. Suppose that $s$ and
$a^{(1)}, \ldots, a^{(m)}$ have $p$-ary representations $s =
\sum_{i=0}^{n-1}s_ip^i$ and $ a^{(j)}= \sum_{i=0}^{n-1}a^{(j)}_ip^i$
for $j=1, \ldots, m$, where the $p$-ary digits $s_i$ and $a^{(j)}_i$
are integers in $\{0,1,\ldots, p-1\}$. Then there exists a {\em
unique\/} integer sequence $c=c_{-1}, c_0, \ldots, c_{n-1}$ with
$c_{-1}=c_{n-1}$ such that \beql{AWCmod} pc_i + s_i = c_{i-1} +
\sum_{j=1}^m t_ja_i^{(j)} \quad (0\leq i\leq n-1).\eeql Moreover, if
we define
\[ t_+=\sum_{\stackrel{j=1}{t_j>0}}^{m}t_j,\qquad
t_-=\sum_{\stackrel{j=1}{t_j<0}}^{m}t_j,\] then $t_--1\leq c_i \leq
t_+$, and furthermore \beql{AWCmodbound} t_-\leq c_i \leq t_+-1
\eeql for $i=0, \ldots, n-1$ provided that $a^{(j)}\not\equiv0\bmod
p^n-1$ for some $j=1, \ldots, m$. As a consequence, the value
$w(c)=c_0+\cdots+c_{n-1}$, the weight $w(s)=s_0+\cdots +s_{n-1}$ of
$s$, and the weights $w(a^{(j)})= a^{(j)}_0+\cdots +a^{(j)}_{n-1}$
of the $a^{(j)}$ satisfy \beql{AWCmodc} (p-1)w(c)=\sum_{j=1}^m
t_jw(a^{(j)}) - w(s). \eeql
\end{teor}

We will usually refer to the $s_i$ and $c_i$ as the ($p$-ary) {\em
digits\/} and {\em carries\/} for the computation modulo $p^n-1$ of
the number $s$. We emphasize that the non-obvious part of
Theorem~\ref{TAWCmod} is the existence of a carry sequence with
$c_{n-1}=c_{-1}$: otherwise (\ref{AWCmod}) represents the ordinary
$p$-ary add-with-carry algorithm. To stress the periodic nature of
this {\em modular\/} add-with-carry algorithm, we will often
consider all indices modulo $n$.

There are various ways to prove this theorem. Here we will derive it from the
following simple technical lemma.
\begin{lem}\label{fawc} Let $r_0, r_1, \ldots, r_{n-1}$ be an integer
sequence.
For all $j$, write
\[r(j)=\sum_{i=0}^{n-1}r_{i+j}p^i,\]
where the indices are to be interpreted modulo $n$.
Then there exists an integer sequence
$c=c_{-1}, c_0, c_1, \ldots, c_{n-1}$ with $c_{-1}=c_{n-1}$ such that
\beql{fawc1} pc_i=r_i+c_{i-1}\eeql
for $i=0, \ldots, n-1$ if and only if $r(0)\equiv 0 \bmod p^n-1$.
In that case, we have
\beql{fawc2}c_{j-1}=r(j)/(p^n-1)\eeql
for $j=0, \ldots, n-1$; in particular, the solution is {\em unique\/}.
Moreover, the ``weights''
$w(r)=r_0+\cdots+r_{n-1}$ and $w(c)=c_0+\cdots+c_{n-1}$
of $r$ and $c$ satisfy
\beql{fawc3}(p-1)w(c) = w(r).\eeql
\end{lem}
\bpf Suppose that (\ref{fawc1}) holds for $i=0, \ldots, n-1$,
with $c_{-1}=c_{n-1}$. Write
\[c(j)=\sum_{i=0}^{n-1} c_{i+j}p^i,\]
where the indices of $c$ are to be interpreted modulo $n$.
Then
\begin{eqnarray*}
pc(j) &=& r(j) + \sum_{i=0}^{n-1} c_{i-1 +j}p^i \\
         &=& r(j) + pc(j) +c_{j-1} -p^nc_{j+n-1} \\
         &=& r(j) + pc(j) - c_{j-1}(p^n-1),
\end{eqnarray*}
hence $r(j)= c_{j-1} (p^n-1)$ for all $j$. So $r(j) \equiv0\bmod p^n-1$ and
$c_{j-1}= r(j)/ (p^n-1)$ for
all $j$; in particular, we have that $r(0)\equiv0 \bmod p^n-1$.

Conversely, suppose that $r(0)\equiv0 \bmod p^n-1$. Then obviously
$p^jr(j) \equiv r(0) \equiv0 \bmod p^n-1$,
so that $r(j)\equiv0\bmod p^n-1$ for all $j$.
Hence the sequence $c_{-1}, c_0, \ldots, c_{n-1}$ defined by (\ref{fawc2})
is an integer sequence, with $c_{-1}=c_{n-1}$ by definition,
and it is easily verified that this sequence indeed satisfies (\ref{fawc1})
for $i=0, \ldots, n-1$.
Finally, the equation (\ref{fawc3}) follows directly from (\ref{fawc1}) by
summing (\ref{fawc1}) for $i=0, \ldots, n-1$.
\epf

\begin{rem}\label{Rpol} If we associate polynomials $r(x)=r_0+r_1x+\cdots
+r_{n-1}x^{n-1}$ and $c(x)=c_0+c_1x+\cdots +c_{n-1}x^{n-1}$ with the
sequences $r$ and $c$, then (\ref{fawc1}) can be read as
\[ r(x)+(x-p)c(x)\equiv 0 \bmod x^n-1.\]
If $p$ is not a zero of $x^n-1$, that is, if $p^n\neq1$, then for each $r$
there is a unique solution $c$. Since
$\gc(x)=(p^{n-1}+p^{n-2}x+\cdots+px^{n-2}+x^{n-1})/(p^n-1)$ satisfies
\[(p-x) \gc(x)\equiv 1 \bmod x^n-1, \]
the solution $c$ is given by
\[ c(x) =\sum_{i=0}^{n-1}r_i \gc(x) x^i.\]
This approach provides an alternative proof of the lemma.
\end{rem}

\bpft{\ref{TAWCmod}}
Define $r_i=-s_i+\sum_{j=1}^kt_ja^{(j)}_i$ for $i=0, \ldots, n-1$.
Since
\[r(0)=-s+\sum_{j=1}^k t_ja^{(j)}\equiv 0 \bmod p^n-1,\]
the existence and uniqueness of the carry sequence $c_{-1}, \ldots, c_{n-1}$
satisfying
(\ref{AWCmod}), as well as the relation (\ref{AWCmodc}),
follows from Lemma~\ref{fawc}.
To obtain the bounds on the carries $c_i$, simply note that
\[(p^n-1)t_- -(p^n-1) \leq r(j) \leq (p^n-1)t_+\]
holds for all $j$. Moreover, since all $t_j$ are assumed to be nonzero,
if equality holds in either of these bounds
then $s$ and each $a^{(j)}$ is equal to 0 modulo $p^n-1$.
\epft

%
%
\section{The weight inequalities}
We will now use Theorem~\ref{TAWCmod} for a local analysis of the
weight inequalities (\ref{wtinequ1}) and (\ref{wtinequ2}). We begin
by analyzing the ternary representations of the two constants
\beql{Euv} u = (3^{2k}-1)/4, \qquad v = 3(3^{2k}-1)/4 \eeql
occurring in (\ref{wtinequ1}) and (\ref{wtinequ2}). Write
\[z=(3^{2k}-1)/8.\]
\begin{lem}\label{Lconst}
The numbers $z,u,v$ are all integers, with $0\leq z,u,v< 3^{2k}-1$,
and $u=2z$, $v=3u$, and $3u\equiv -u\equiv v \bmod 3^{2k}-1$.
Moreover, if $z=z_{2k-1}\cdots z_0$, $u=u_{2k-1}\cdots u_0$, and
$v=v_{2k-1}\cdots v_0$ are the (unique) $3$-ary representations of
$z$, $u$, and $v$, respectively, then
\[ z_{2i}=1,\qquad z_{2i-1}=0, \qquad u_i = 2z_i, \qquad v_i=u_{i-1}.\]
In particular, $w(u)=w(v)=2w(z)=2k$. Finally, if $k$ is odd, then
$3^k u\equiv v$ and $3^kv\equiv u$ modulo $3^{2k}-1$.
\end{lem}
\bpf By definition, $z=(3^{2k}-1)/8=1+3^2 + \cdots+3^{2k-2}$, so $z$
is integer, and hence $u=2z$ and $v=3u$ are also integers. Since
also $0<z<u<v<3^{2k}-1$, their $3$-ary representations are unique
and as described in the lemma. This immediately implies that
$w(u)=w(v)=2w(z)=2k$. Next note that $4u = u+v = 3^{2k}-1 \equiv 0
\bmod 3^{2k}-1$; hence $3u \equiv -u \equiv v$ and $3v \equiv
-v\equiv u$ modulo $3^{2k}-1$. As a consequence, if $k$ is odd, then
$3^ku\equiv 3u\equiv -u \equiv v\bmod 3^{2k}-1$ and $3^kv\equiv
u\bmod 3^{2k}-1$. \epf

Now let $b$ be any integer. Suppose that $k$ is odd. Then $3^ku=v$
by Lemma~\ref{Lconst}. Hence if  $s$ satisfies
\[ s\equiv -(3^k+1)b -u \bmod 3^{2k}-1,\]
then $t=3^ks$ satisfies
\[ t \equiv -(3^k+1)b-v \bmod 3^{2k}-1.\]
Now the ternary representation of $t$ is just a cyclic shift of that
of $s$, hence $w(s)=w(t)$ (this is also correct if $s\equiv t\equiv
0\bmod 3^{2k}-1$), and so the two weight inequalities
(\ref{wtinequ1}) and (\ref{wtinequ2}) are in fact equivalent.

We want to prove these two weight inequalities by analyzing the
contribution to the weights from individual ternary digits of $b$.
To this end, we will apply Theorem~\ref{TAWCmod} to the addition
$s\equiv -b-a +v \bmod 3^{2k}-1$, where $a=3^kb$, with the aim to
prove the weight inequality $w(b)+w(s)\geq 2k$. (For technical
reasons, we will prefer this form of the addition, which has the
same outcome $s$ as the earlier one since $-u\equiv v \bmod
3^{2k}-1$). However, the $i^{\rm th}$ digit $a_i=b_{i+k}$ (indices
modulo $2k$) of $a$ and the $i^{\rm th}$ digit $b_i$ of $b$ are
entirely unrelated, so a straightforward local analysis is doomed to
fail. This is a standard problem when investigating $p$-ary weight
inequalities, see for example several cases in \cite{hx}.
Fortunately, there is a standard solution. First, use the relations
between the numbers involved to derive other, equivalent weight
inequalities, typically by multiplying the relevant addition by a
suitable power of $p$. Then forget the relation between the numbers,
but analyze the resulting weight inequalities (in fact, their
addition) {\em simultaneously\/}. This approach in general results
in a {\em generalisation\/} of the original weight inequality.

For the case at hand, the generalization suggested by this approach
turns out to be the following.
\begin{teor}\label{Tgenwi}
Let $k$ be any positive integer, and let $u$ and $v$ be defined as in
(\ref{Euv}). For any integers $a$ and $b$, if $s$ and $t$ satisfy
\[ s \equiv -a-b+v, \qquad t \equiv -a-b+u \bmod 3^{2k}-1, \]
then $w(a)+w(b)+w(s)+w(t)\geq 4k$.
\end{teor}
We first show that the original weight inequalities (\ref{wtinequ1}) and
(\ref{wtinequ2}) indeed follow from this result.
\begin{cor}\label{Cwi}
Let $b$ be an integer, and let $k$ be odd. Then
\[ w(b)+w(-(3^k+1)b-(3^{2k}-1)/4)=w(b)+w(-(3^k+1)b-3(3^{2k}-1)/4)\geq 2k.\]
\end{cor}
\bpf
Apply Theorem~\ref{Tgenwi} with $a=3^kb$. Then $w(a)=w(b)$. Also,
since $k$ is odd, $3^ku\equiv v\bmod 3^{2k}-1$ by Lemma~\ref{Lconst}; hence
\[ 3^ks \equiv -3^k(3^k+1)b +3^k v \equiv -(3^k+1)b+u \equiv t \bmod
3^{2k}-1,\]
so that $w(t)=w(s)$. From the theorem, we now conclude that
$w(b)+w(s)=w(b)+w(t)\geq 2k$, as claimed.
\epf

\bpft{\ref{Tgenwi}} Let $s$ and $t$ be defined as in the theorem,
and write $n=2k$. Assume that $a$, $b$, $u$, $v$, $s$, and $t$ have
ternary digits $a_i$, $b_i$, $u_i$, $v_i$, $s_i$, and $t_i$, for
$i=0, \ldots, n-1$. According to Lemma~\ref{Lconst}, we have that
$u_i=v_{i-1}$ (indices modulo $n$) and $v_{2i}=0$, $v_{2i-1}=2$, for
all $i$. Now apply Theorem~\ref{TAWCmod} to the defining additions
for $s$ and $t$. In both cases, $t_+=1$ and $t_-=2$, hence there are
carry sequences $c_0, \ldots, c_{n-1}$ and $d_0, \ldots, d_{n-1}$
with $c_{-1}=c_{n-1}$, $d_{-1}=d_{n-1}$, and $-2\leq c_i, d_i \leq
0$ for all $i$ such that
\beqal{acforst} 3 c_i + s_i &=& -a_i -b_i +v_i+c_{i-1}\\
      3 d_i + t_i &=& -a_i -b_i +v_{i-1}+d_{i-1}
\eeqal for all $i=0, \ldots, n-1$. Moreover,
$2w(c)+w(s)=2w(d)+w(t)=-w(a)-w(b)+w(v)=-w(a)-w(b)+2k$. Using these
relations, we see that the weight inequality in the theorem is
equivalent to \beql{Ewalt} w(a)+w(b)+2w(c)+2w(d)\leq 0.\eeql

In order to analyze the contribution of the individual ternary
digits $a_i$, $b_i$ to the sum of the weights in the left hand side
of (\ref{Ewalt}), we construct the following labeled directed graph
$G$.

The graph $G$ will have a vertex
$ (a',b',c',d',v')$ whenever $a',b'\in \{0,1,2\}$, $c',d'\in \{-2,-1,0\}$, and
$v'\in\{0,2\}$, and a weighted directed arc
\[ (a',b',c',d',v')\ \ \stackrel{a'+b'+2c''+2d''}{\longrightarrow}\ \
(a'',b'',c'',d'',v'')\]
whenever $v''=2-v'$,
\[ s'= -a'-b'+v'+c'-3c'' \in \{0,1,2\}, \qquad
t'=-a'-b'+v''+d'-3d'' \in \{0,1,2\}.\]
Note that, according to these definitions, whenever
(\ref{acforst}) holds there is an arc
\[ (a_i,b_i,c_{i-1},d_{i-1},v_{i-1})\ \
\stackrel{a_i+b_i+2c_i+2d_i}{\longrightarrow}\ \
(a_{i+1},b_{i+1},c_i,d_i,v_i)\] in the graph. Moreover, since
$c_{-1}=c_{n-1}$ and $d_{-1}=d_{n-1}$, there is a one-to-one
correspondence between sets of relations (\ref{acforst}) for $i=0,
\ldots, n-1$ with corresponding sum of weights
$w=w(a)+w(b)+2w(c)+2w(d)$ and directed walks of length $n$ in the
graph for which the sum of the weights of the arcs equals $w$.

So we are done if we can show that the weight of each directed walk
in the graph $G$ is non-positive. We will show that in fact {\em
all\/} arc-weights are non-positive. To this end, consider an arc
\[ (a',b',c',d',v')\ \ \stackrel{a'+b'+2c''+2d''}{\longrightarrow}\ \
(a'',b'',c'',d'',v''),\]
where $a',b'\in \{0,1,2\}$, $c',d', c'', d''\in \{-2,-1,0\}$,
$v'\in\{0,2\}$, $v''=2-v'$, and
\[ s'= v' -(a'+b')-3c''+c' \in \{0,1,2\}, \qquad
t'=v''-(a'+b')-3d''+d' \in \{0,1,2\}.\] Since $s',t'\geq0$ and
$c',d'\leq0$, we conclude that \beql{Earcw} 3c''\leq v'-(a'+b'),
\qquad 3d''\leq v''-(a'+b'). \eeql Now consider the arc weight
$w=(a'+b')+2(c''+d'')$. We have that $0\leq a'+b'\leq 4$. Obviously,
if $a'+b'=0$, then $w\leq0$. If $1\leq a'+b'\leq2$, then since
$v'=0$ or $v''=0$, at least one of $c''$ or $d''$ is negative, and
$w\leq0$ again. Finally, if $3 \leq a'+b'\leq4$, then both
$c'',d''\leq -1$, and again $w\leq0$. \epft

A close inspection of this proof reveals that in the same way, we can prove the
following somewhat stronger result.
\begin{teor}\label{Tgengenwi}
Let $n$ be any positive integer. Let $u$ and $v$ be numbers with $3$-ary
representation $u=u_{n-1}\cdots u_0$ and $v=v_{n-1}\cdots v_0$, respectively,
for which $u_i=0$ or $v_i=0$ holds for each $i=0, \ldots, n-1$.
If $a$, $b$ are integers, and $s$ and $t$ satisfy
\[ s \equiv -a-b+v, \qquad t \equiv -a-b+u \bmod 3^n-1, \]
then $w(a)+w(b)+w(s)+w(t)\geq w(u)+w(v)$.
\end{teor}
\bibliographystyle{amsalpha}

\end{document}